# Spline Based Series for Sine and Arbitrarily Accurate Bounds for Sine, Cosine and Sine Integral

**Roy M. Howard**

School of Electrical Engineering, Computing and Mathematical Sciences, Curtin University, GPO Box U1987, Perth, 6845, Australia.

email: r.howard@curtin.edu.au

**17 June 2020**

**Abstract**

Based on two point spline approximations of arbitrary order, a series of functions that define lower bounds for $\sin(x)$ and $\sin(x)/x$, over the interval $[0, \pi/2]$, with increasingly low relative errors and smaller relative errors than published results, are defined. Second, fourth and eighth order approximations have, respectively, maximum relative errors over the interval $[0, \pi/2]$ of $3.31 \times 10^{-4}$, $2.48 \times 10^{-8}$ and $2.02 \times 10^{-18}$. New series for the sine function, which have significantly better convergence that a Taylor series over the interval $[0, \pi/2]$, are proposed. Applications include functions that are upper bounds for the sine function, upper and lower bounds for the cosine function and lower bounds for the sine integral function. These bounded functions can be made arbitrarily accurate.

# 1   Introduction

Research to establish upper and lower bounds for $\sin(x)/x$ over the interval $[0, \pi/2]$ has a long history (e.g. Sándor & Bencze 2005, Mortici 2011) and Qi et al. 2009 provides a detailed review. Fundamental results include Jordan's inequality

$$\frac{2}{\pi} \leq \frac{\sin(x)}{x} < 1 \qquad x \in [0, \pi/2] \tag{1}$$

the Cusa-Huygens inequality

$$\frac{\sin(x)}{x} < \frac{2 + \cos(x)}{3} \qquad 0 < x < \frac{\pi}{2} \tag{2}$$

and the Redheffer inequality (Redheffer 1969):

$$\frac{\pi^2 - x^2}{\pi^2 + x^2} \leq \frac{\sin(x)}{x} \qquad x \in [0, \pi] \tag{3}$$

Similar, and related inequalities, which are indicative of published bounds, are detailed in Table 1.1. The graphs of the relative errors in the bounds defined by the first, second, fourth, fifth, eighth and tenth inequalities stated in this table are shown in Figure 1. The relative error of the other bounds specified in Table 1.1 lie within the limits of the bounds shown in Figure 1. The relative errors associated with these bounds typical have maximum values, over the interval $[0, \pi/2]$, in the range of $10^{-4}$ to $0.2$. Some of the bound are 'sharp', i.e. exact, at the points zero and $\pi/2$. The majority of the bounds have relative errors that are monotonically increasing and have maximum values at $\pi/2$.

Applications of such bounds include defining bounds to the Riemann zeta function $\zeta(3)$, e.g. Luo et al., 2003, and bounds on some trigonometric ratios, e.g. Wu, 2004.

**Table 1.1**   Published upper and lower bounds for $\sin(x)/x$, $0 < x < \pi/2$.

| # | Reference | Lower bound | Upper bound |
|---|---|---|---|
| 1 | Lower: Qi, 1999, eqn. 26. Upper: Cusa-Huygens | $\dfrac{1 + \cos(x)}{2}$ | $\dfrac{2 + \cos(x)}{3}$ |
| 2 | Lower: Mitrinović, 1965 Upper: Cusa-Huygens | $\cos(x)^{1/3}$ | $\dfrac{2 + \cos(x)}{3}$ |
| 3 | Bhayo, 2017, eqn. 1.6 | $\dfrac{\cos(x) + \alpha - 1}{\alpha}, \alpha = \dfrac{\pi}{\pi - 2}$ | $\dfrac{\cos(x) + \beta - 1}{\beta}, \beta = 3$ |
| 4 | Bercu, 2016, Lemma 2.1 | $\dfrac{1 - 7x^2/60}{1 + x^2/20}$ | $\dfrac{1 - x^2/7 + 11x^4/2520}{1 + x^2/42}$ |
| 5 | Hua, 2016, eqn. 3.1 | $2 + \dfrac{23x^3 \sin(x)}{720} - \dfrac{\tan(x/2)^2}{(x/2)^2}$ | $2 + \dfrac{k_0 x^3 \sin(x)}{\pi^5} - \dfrac{\tan(x/2)^2}{(x/2)^2}$ $k_0 = 128 - 16\pi^2 + 16\pi$ |
| 6 | Yang, 2016, eqn. 2.14 | $\left[\dfrac{2}{\pi}\right]^{4x^2/\pi^2}$ | $\exp\left[\dfrac{-x^2}{6}\right]$ |
| 7 | Yang, 2014, eqn. 54 | $\cos[p_o x]^{1/p_o}, p_o = 0.3473$ | $\cos\left[\dfrac{x}{3}\right]^3$ |
| 8 | Chen, 2015, eqn 4.7 | $\dfrac{p + 6\cos(x)}{14 + \cos(x)}, p = \dfrac{28}{\pi}$ | $\dfrac{9 + 6\cos(x)}{14 + \cos(x)}$ |
| 9 | Chen, 2015, eqn 4.9 | $\left[\dfrac{9 + 6\cos(x)}{14 + \cos(x)}\right]^r, r = \dfrac{\ln(\pi/2)}{\ln(14/9)}$ | $\dfrac{9 + 6\cos(x)}{14 + \cos(x)}$ |
| 10 | Chen, 2015, eqn 5.8 | $\dfrac{2 + \cos(x) - k_0 x^2}{3 - k_0 x^2}, k_0 = \dfrac{8\pi - 24}{\pi^3 - \pi^2}$ | $\dfrac{2 + \cos(x) - k_1 x^2}{3 - k_1 x^2}, k_1 = \dfrac{1}{10}$ |

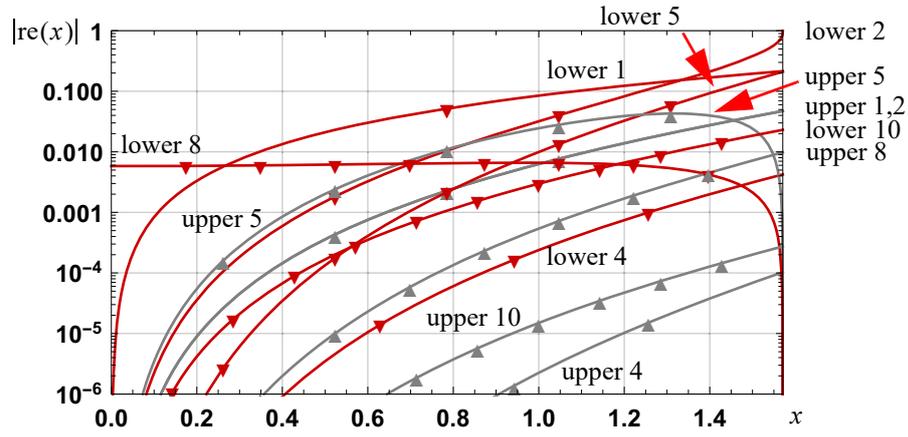

**Figure 1** Graph of the magnitude of the relative errors in the upper and lower bounds to $\sin(x)/x$ as defined by the first, second, fourth, fifth, eighth and tenth order functions defined in Table 1.1.

Naturally, it is of interest to determine bounds with lower relative errors. Zhu (2008a, Theorem 5) proposed the following general form for bounds of $\sin(x)/x$, $0 < x < \pi/2$:

$$\sum_{k=0}^{n} \alpha_k (\pi^2 - 4x^2)^k + \alpha_{n+1}(\pi^2 - 4x^2)^{n+1} \leq \frac{\sin(x)}{x} \leq \sum_{k=0}^{n} \alpha_k (\pi^2 - 4x^2)^k + \left[1 - \sum_{k=0}^{n} a_k \pi^{2k}\right] \cdot \frac{(\pi^2 - 4x^2)^{n+1}}{\pi^{2n+2}} \quad (4)$$

$$\alpha_k = \frac{2k-1}{2k\pi^2} \cdot \alpha_{k-1} - \frac{1}{16(k-1)k\pi^2} \cdot \alpha_{k-2} \qquad \alpha_0 = \frac{2}{\pi} \qquad \alpha_1 = \frac{1}{\pi^3} \qquad n \in \{0, 1, 2, \ldots\}$$

and, as indicated in Figure 2, the relative error in these bounds can be made arbitrarily small by utilizing increasingly high orders. This result was generalized in Zhu 2008b. These bounds build on the results of Li, 2006 (Theorem 2.1) and an alternative, but equivalent form, was published by Nui et al. (2008). Nui et al. 2010 provides a generalization and showed that the equation forms specified by Nui and Zhu are equivalent (Proposition 3). These bounds, for orders zero to two, are:

$$\frac{2}{\pi} + \frac{\pi^2 - 4x^2}{\pi^3} < \frac{\sin(x)}{x} < \frac{2}{\pi} + \frac{1 - 2/\pi}{\pi^2} \cdot (\pi^2 - 4x^2) \quad (5)$$

$$\frac{2}{\pi} + \frac{\pi^2 - 4x^2}{\pi^3} + \frac{12 - \pi^2}{16\pi^5} \cdot (\pi^2 - 4x^2)^2 < \frac{\sin(x)}{x} < \frac{2}{\pi} + \frac{\pi^2 - 4x^2}{\pi^3} + \frac{1 - 3/\pi}{\pi^4} \cdot (\pi^2 - 4x^2)^2 \quad (6)$$

$$\frac{2}{\pi} + \frac{\pi^2 - 4x^2}{\pi^3} + \frac{12 - \pi^2}{16\pi^5} \cdot (\pi^2 - 4x^2)^2 + \frac{10 - \pi^2}{16\pi^7} \cdot (\pi^2 - 4x^2)^3 < \frac{\sin(x)}{x} <$$

$$\frac{2}{\pi} + \frac{\pi^2 - 4x^2}{\pi^3} + \frac{12 - \pi^2}{16\pi^5} \cdot (\pi^2 - 4x^2)^2 + \frac{1}{\pi^6} \cdot \left[1 - \frac{3}{\pi} - \frac{12 - \pi^2}{16\pi}\right] \cdot (\pi^2 - 4x^2)^3 \quad (7)$$

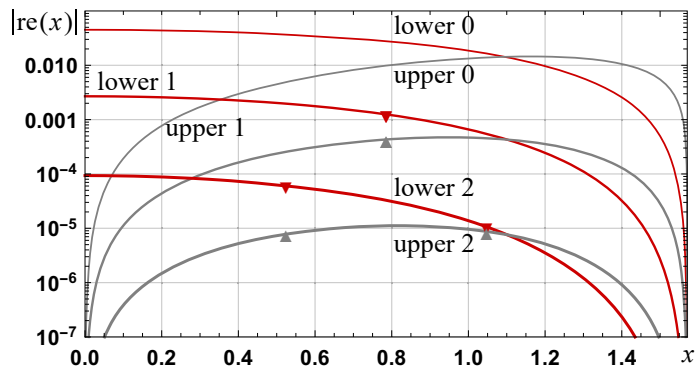

**Figure 2** Relative error in upper and lower bounds to $\sin(x)/x$, of orders zero to two, as defined by Equation 5 to Equation 7.

This paper utilizes the two point spline based function approximation, detailed in Howard 2019, to specify a series of upper and lower bounds for $\sin(x)$ and $\sin(x)/x$, over the interval $[0, \pi/2]$, with increasingly low relative errors. Analysis associated with the bounds leads to new series for the sine function which have significantly better convergence, over the interval $[0, \pi/2]$, than a Taylor series and the published series proposed by Zhu 2008a and Nui 2008. Applications include functions that define upper and lower bounds for the cosine function and lower bounds for the sine integral function. These are can be made arbitrarily accurate.

Section 2 details miscellaneous results that underpin subsequent analysis. Section 3 details a spline based series for the sine function and in Section 4 it is shown that these series are lower bounds, of increasing accuracy. Section 5 details applications including new series for the sine function, upper bounds for the sine function, bounds for the cosine function and a lower bound for the sine integral function. Conclusions are detailed in section 6.

## 1.1   Notation

For a function $f$ defined over the interval $[\alpha, \beta]$, an approximating function $f_A$ has a relative error, at a point $x_1$, defined according to

$$\text{re}(x_1) = 1 - \frac{f_A(x_1)}{f(x_1)} \tag{8}$$

The relative error bound for the approximating function $f_A$, over the interval $[\alpha, \beta]$, is defined according to

$$\text{re}_B = \max\{\text{re}(x_1): x_1 \in [\alpha, \beta]\} \tag{9}$$

The notation $f^{(k)}(x) = \frac{d^k}{dt^k} f(x)$ is used.

Mathematica has been used to facilitate analysis and to obtain numerical results. In general, relative error bound results for approximations have been obtained by sampling the appropriate interval at 1000 equally spaced points.

## 2   Miscellaneous Results

The following results underpin later analysis and discussion:

### 2.1   Taylor Series

The default reference for bounds for the sine function, and hence $\sin(x)/x$, is the Taylor series defined according to

$$\sin(x) = x - \frac{x^3}{3!} + \frac{x^5}{5!} - \frac{x^7}{7!} + \frac{(-1)^{\lfloor k/2 \rfloor} x^k}{k!} + \ldots \quad k \in \{1, 3, 5, 7, \ldots\} \tag{10}$$

A Taylor series yields alternating lower and upper bounds over the interval $[0, \pi/2]$ as the order is increased (e.g. Qi, 2009, Section 1.6). The bounds are increasingly accurate as can be seen in the results shown, along with other results, in Figure 3. The relative error bounds, over the interval $[0, \pi/2]$, for various orders of Taylor series approximations, are detailed in Table 2.1. The rate of convergence of the Taylor series serves as a reference.

Table 2.1   Relative error bounds for Taylor series approximations to $\sin(x)$ over the interval $[0, \pi/2]$.

| Order of approx. | Relative error bound |
|---|---|
| 1: upper bound | 0.571 |
| 3: lower bound | $7.52 \times 10^{-2}$ |
| 5: upper bound | $4.52 \times 10^{-3}$ |

Table 2.1   Relative error bounds for Taylor series approximations to $\sin(x)$ over the interval $[0, \pi/2]$.

| Order of approx. | Relative error bound |
|---|---|
| 7: lower bound | $1.57 \times 10^{-4}$ |
| 9: upper bound | $3.54 \times 10^{-6}$ |
| 13: lower bound | $6.63 \times 10^{-10}$ |
| 17: lower bound | $4.35 \times 10^{-14}$ |
| 33: upper bound | $7.07 \times 10^{-34}$ |

## 2.2   Generating an Upper Bound from a Lower Bound

The following lemma details when upper bounds to a set function $f$ can be generated from a sequence of converging lower bounds $f_1^L, f_2^L, \ldots$:

**Lemma 1**   **Upper Bound from a Lower Bound**

Given a sequence of convergent lower bounds $f_1^L, f_2^L, \ldots$ to a function $f$ over an interval $[\alpha, \beta]$, a sequence of upper bounds can be defined over the interval $[\alpha, \beta]$ according to

$$f_k^U(x) = 2f_k^L(x) - f_{k-1}^L(x) \qquad x \in [\alpha, \beta], k \in \{2, 3, \ldots\} \tag{11}$$

provided $\varepsilon_{k-1}^L(x) > 2\varepsilon_k^L(x)$, $x \in [\alpha, \beta]$ where

$$\varepsilon_k^L(x) = f(x) - f_k^L(x) \qquad x \in [\alpha, \beta] \tag{12}$$

**Proof**

Consider $f_k^U(x) = 2f_k^L(x) - f_{k-1}^L(x)$ and the error

$$\begin{aligned}\varepsilon_k^U(x) &= f_k^U(x) - f(x) = 2f_k^L(x) - f_{k-1}^L(x) - f(x) \\ &= [2f_k^L(x) - 2f(x)] + [f(x) - f_{k-1}^L(x)] = -2\varepsilon_k^L(x) + \varepsilon_{k-1}^L(x)\end{aligned} \tag{13}$$

When $\varepsilon_{k-1}^L(x) > 2\varepsilon_k^L(x)$ it follows that $f_k^U(x) > 0$ as required.

## 2.3   Relationship Between Sine and Cosine Bounds

**Lemma 2**   **Relationship Between Sine and Cosine Bounds**

The bounds $f_L(x) < \sin(x) < f_U(x)$, $x \in [0, \pi/2]$ for the sine function imply the following bounds for the cosine function:

$$f_L\left[\frac{\pi}{2} - y\right] < \cos(y) < f_U\left[\frac{\pi}{2} - y\right] \qquad y \in \left[0, \frac{\pi}{2}\right] \tag{14}$$

**Proof**

As $\sin(\pi/2 - x) = \cos(x)$ and $\sin(x) = \cos(\pi/2 - x)$, for $x \in [0, \pi/2]$, it follows that

$$f_L(x) < \sin(x) < f_U(x) \quad \Rightarrow \quad f_L(x) < \cos\left[\frac{\pi}{2} - x\right] < f_U(x) \qquad x \in \left[0, \frac{\pi}{2}\right]$$

$$\Rightarrow \quad f_L\left[\frac{\pi}{2} - y\right] < \cos(y) < f_U\left[\frac{\pi}{2} - y\right] \qquad y = \frac{\pi}{2} - x, x = \frac{\pi}{2} - y \qquad y \in \left[0, \frac{\pi}{2}\right]$$

(15)

### 2.4   Function Scaling

When considering analysis of the sine function the basic choices are between considering $\sin(x)$ over the interval $[0, \pi/2]$ or $\sin[\pi x/2]$ over the interval $[0, 1]$. When establishing approximations to the sine function, the former is preferred whilst the latter facilitates establishing convergence of approximations. The appropriate relationships between the two forms are detailed below:

**Lemma 3      Scaling Relationships**

With approximating functions $f_0$ and $f_1$, and error functions $\varepsilon_0$ and $\varepsilon_0$, defined according to

$$y_0(x) = \sin(x) = f_0(x) + \varepsilon_0(x) \qquad x \in \left[0, \frac{\pi}{2}\right]$$

$$y_1(x_1) = \sin\left[\frac{\pi}{2} \cdot x_1\right] = f_1(x_1) + \varepsilon_1(x_1) \qquad x_1 \in [0, 1]$$

(16)

the linear transformations $x = \pi x_1/2$ and $x_1 = 2x/\pi$, for the case of $f_0(x) = f_1[2x/\pi]$, $f_1(x_1) = f_0[\pi x_1/2]$, are

$$y_0(x) = y_1\left[\frac{2}{\pi} \cdot x\right] \qquad y_1(x_1) = y_0\left[\frac{\pi}{2} \cdot x_1\right]$$

$$\varepsilon_0(x) = \varepsilon_1\left[\frac{2}{\pi} \cdot x\right] \qquad \varepsilon_1(x_1) = \varepsilon_0\left[\frac{\pi}{2} \cdot x_1\right]$$

(17)

The relative errors in the approximations $\sin(x) \approx f_0(x_1)$ and $\sin(\pi x_1/2) \approx f_1(x_1)$, respectively, are:

$$re_0(x) = 1 - \frac{f_0(x)}{\sin(x)} \qquad re_1(x_1) = 1 - \frac{f_1(x_1)}{\sin(\pi x_1/2)}$$

(18)

and the relationships between the two relative errors are

$$re_0(x) = re_1(2x/\pi) \qquad re_1(x_1) = re_0(\pi x_1/2)$$

(19)

**Proof**

Consider, for example:

$$re_0(\pi x_1/2) = 1 - \frac{f_0[\pi x_1/2]}{\sin[\pi x_1/2]} = 1 - \frac{f_1(x_1)}{\sin[\pi x_1/2]} = re_1(x_1)$$

(20)

#### 2.4.1    Note

It is clear that establishing the relative error for an approximation to $\sin(\pi x_1/2)$ for $x_1 \in [0, 1]$ is equivalent to establishing the relative error for the scaled approximation to $\sin(x)$ for $x \in [0, \pi/2]$.

## 3    Spline Based Series for Sine

A $n$th order, two point, spline approximation to a function $f$, which is at least a $n$th order differentiable over the interval $[\alpha, \beta]$, has been detailed in Howard 2019 and is

$$f_n(x) = \frac{(\beta-x)^{n+1}}{(\beta-\alpha)^{n+1}} \cdot \sum_{k=0}^{n} \frac{(x-\alpha)^k}{k!} \cdot f^{(k)}(\alpha) \cdot \left[\sum_{i=0}^{n-k} \frac{(n+i)!}{i! \cdot n!} \cdot \frac{(x-\alpha)^i}{(\beta-\alpha)^i}\right] +$$

$$\frac{(x-\alpha)^{n+1}}{(\beta-\alpha)^{n+1}} \cdot \sum_{k=0}^{n} \frac{(-1)^k(\beta-x)^k}{k!} \cdot f^{(k)}(\beta) \cdot \left[\sum_{i=0}^{n-k} \frac{(n+i)!}{i! \cdot n!} \cdot \frac{(\beta-x)^i}{(\beta-\alpha)^i}\right] \qquad x \in [\alpha, \beta] \quad (21)$$

This general result can be utilized to define a sequence of functions that approximate, with increasing accuracy, the sine function:

### Theorem 3.1   Spline Series for Sine

A $n$th order spline series approximation to $\sin(x)$, for the interval $[0, \pi/2]$, and based on the points $0$ and $\pi/2$, is

$$f_n(x) = \left[1 - \frac{2x}{\pi}\right]^{n+1} \sum_{k=0}^{n} \frac{1}{k!} \cdot \sin\left[\frac{k\pi}{2}\right] \cdot x^k \cdot \sum_{i=0}^{n-k} \frac{(n+i)!}{i! \cdot n!} \cdot \left[\frac{2x}{\pi}\right]^i +$$

$$\left[\frac{2x}{\pi}\right]^{n+1} \sum_{k=0}^{n} \frac{(-1)^k}{k!} \cdot \left[\frac{\pi}{2}\right]^k \cdot \sin\left[\frac{(k+1)\pi}{2}\right] \cdot \left[1 - \frac{2x}{\pi}\right]^k \cdot \sum_{i=0}^{n-k} \frac{(n+i)!}{i! \cdot n!} \cdot \left[1 - \frac{2x}{\pi}\right]^i \qquad (22)$$

$$n \in \{0, 1, 2, \ldots\}$$

The approximations, of order zero to fourth, are:

$$f_0(x) = \frac{2x}{\pi} \qquad (23)$$

$$f_1(x) = x + \frac{12}{\pi^2}\left[1 - \frac{\pi}{3}\right]x^2 - \frac{16}{\pi^3}\left[1 - \frac{\pi}{4}\right]x^3 \qquad (24)$$

$$f_2(x) = x + \frac{80}{\pi^3}\left[1 - \frac{3\pi}{10} - \frac{\pi^2}{80}\right]x^3 - \frac{240}{\pi^4}\left[1 - \frac{4\pi}{15} - \frac{\pi^2}{60}\right]x^4 + \frac{192}{\pi^5}\left[1 - \frac{\pi}{4} - \frac{\pi^2}{48}\right]x^5 \qquad (25)$$

$$f_3(x) = x - \frac{x^3}{6} + \frac{560}{\pi^4}\left[1 - \frac{2\pi}{7} - \frac{\pi^2}{56} + \frac{\pi^3}{420}\right]x^4 - \frac{2688}{\pi^5}\left[1 - \frac{15\pi}{56} - \frac{\pi^2}{48} + \frac{\pi^3}{672}\right]x^5 +$$

$$\frac{4480}{\pi^6}\left[1 - \frac{9\pi}{35} - \frac{13\pi^2}{560} + \frac{\pi^3}{840}\right]x^6 - \frac{2560}{\pi^7}\left[1 - \frac{\pi}{4} - \frac{\pi^2}{40} + \frac{\pi^3}{960}\right]x^7 \qquad (26)$$

$$f_4(x) = x - \frac{x^3}{6} + \frac{4032}{\pi^5}\left[1 - \frac{5\pi}{18} - \frac{\pi^2}{48} + \frac{5\pi^3}{2016} + \frac{\pi^4}{48384}\right]x^5 -$$

$$\frac{26880}{\pi^6}\left[1 - \frac{4\pi}{15} - \frac{11\pi^2}{480} + \frac{\pi^3}{504} + \frac{\pi^4}{40320}\right]x^6 + \frac{69120}{\pi^7}\left[1 - \frac{7\pi}{27} - \frac{53\pi^2}{2160} + \frac{\pi^3}{576} + \frac{\pi^4}{34560}\right]x^7 - \qquad (27)$$

$$\frac{80640}{\pi^8}\left[1 - \frac{16\pi}{63} - \frac{13\pi^2}{504} + \frac{\pi^3}{630} + \frac{\pi^4}{30240}\right]x^8 + \frac{35840}{\pi^9}\left[1 - \frac{\pi}{4} - \frac{3\pi^2}{112} + \frac{\pi^3}{672} + \frac{\pi^4}{26880}\right]x^9$$

Higher order approximations can readily be defined.

### Proof

These results arise from direct application of Equation 21 for the case of $f(x) = \sin(x)$, the interval $[0, \pi/2]$, $\alpha = 0$, $\beta = \pi/2$ and the result

$$f^{(k)}(x) = \sin\left[x + \frac{k\pi}{2}\right] \qquad k \in \{0, 1, 2, \ldots\} \qquad (28)$$

## 3.1  Results

The graphs of the relative errors in the first to fourth order spline based approximations, defined in Theorem 3.1, are shown in Figure 3 and the associated relative error bounds are detailed in Table 3.1. The error functions for the first to fourth order approximations are shown in Figure 4. These results, when compared with Taylor series approximations and the approximations defined by Zhu 2008a (specified in Equation 4 and whose relative errors are shown in Figure 2), demonstrate better convergence. The approximations are exact at the points $0, \pi/2$.

The errors in the approximations, for orders one to four, are shown in Figure 4 and are positive which suggests the approximations can serve as lower bounds for the sine function. The proof of this is detailed below.

**Table 3.1**  Relative error bounds for spline based approximations to $\sin(x)$ over the interval $[0, \pi/2]$.

| Order of approx. | Relative error bound: Spline based. |
|---|---|
| 0 | 0.363 |
| 1 | $1.63 \times 10^{-2}$ |
| 2 | $3.31 \times 10^{-4}$ |
| 3 | $3.62 \times 10^{-6}$ |
| 4 | $2.48 \times 10^{-8}$ |
| 6 | $3.91 \times 10^{-13}$ |
| 8 | $2.02 \times 10^{-18}$ |
| 16 | $9.19 \times 10^{-43}$ |
| 32 | $2.19 \times 10^{-100}$ |

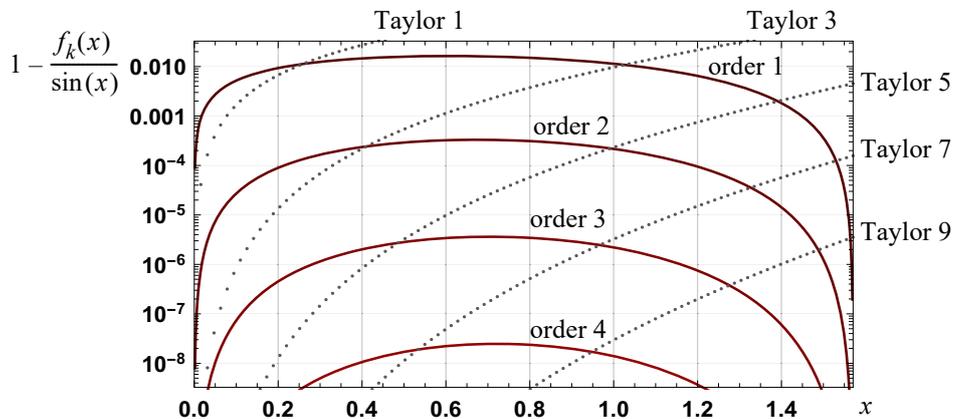

**Figure 3**  Graphs of the relative error in the spline approximations to $\sin(x)$, of orders one to four, as well as Taylor series approximations of orders one, three, five, seven and nine.

## 4  Lower Order Bounds for Sine

To show that the spline approximations for $\sin(x)$, as detailed in Theorem 3.1, are lower bounds, the first and second order approximations are considered.

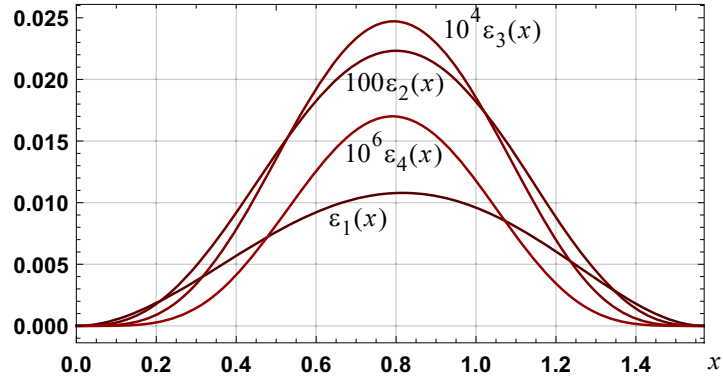

**Figure 4**  Graph of the errors, $\varepsilon_k(x) = \sin(x) - f_k(x)$, in the first to fourth order spline approximations to $\sin(x)$.

## 4.1  Lower Bound Based on First Order Spline Approximation

The proof, that the first order spline approximation, as defined by Equation 24, is a lower bound for $\sin(x)$ over the interval $[0, \pi/2]$, is based on showing that the associated error function is positive over this interval. Consider the scaled error function

$$\varepsilon_1(x_1) = \sin\left[\frac{\pi}{2} \cdot x_1\right] - f_1\left[\frac{\pi}{2} \cdot x_1\right] \qquad x_1 \in [0, 1] \tag{29}$$

which has the property of being positive on the interval $(0, 1)$ and being such that $\varepsilon_1(0) = \varepsilon_1(1) = 0$:

**Theorem 4.1  Error in First Order Spline Approximation to Sine**

Using the notation $t = x_1$, the error function defined by Equation 29, can be written as a summation of positive terms according to

$$\begin{aligned}
\varepsilon_1(t) &= c_2 t^2 (1-t)^2 + c_3 t^3 (1-t)^3 + c_4 t^4 (1-t)^3 + c_5 t^5 (1-t)^2 + c_6 t^6 (1-t)^2 + c_7 t^7 (1-t)^3 + \ldots \\
&= \sum_{k=2}^{\infty} c_k t^k (1-t)^p \qquad p = \begin{cases} 2 & k = 2, 5, 6, 9, 10, 13, 14, 17, 18, \ldots \\ 3 & k = 3, 4, 7, 8, 11, 12, 15, 16, 19, 20 \ldots \end{cases} \\
&= \sum_{k=2}^{\infty} c_k t^k (1-t)^{\left[\frac{5}{2} + \frac{(-1)^{\lfloor (k+1)/2 \rfloor}}{2}\right]} \qquad t \in [0, 1]
\end{aligned} \tag{30}$$

The series converges for $t \in [0, 1]$ and the coefficients can be defined iteratively according to

$$c_0 = -1 \qquad c_1 = -2 + \frac{\pi}{2} = 2c_0 + \frac{\pi}{2} \qquad c_2 = \pi - 3 = 2c_1 - c_0 \tag{31}$$

$$c_3 = -4 + \frac{3\pi}{2} - \frac{\pi^3}{48} = c_2 + c_1 - c_0 - \frac{\pi^3}{2^3 3!} \qquad c_4 = -9 + \frac{7\pi}{2} - \frac{\pi^3}{16} = 3c_3 - c_2 \tag{32}$$

and

$$\begin{aligned}
c_k &= 3c_{k-1} - c_{k-2} & k &\in \{4, 8, 12, 16, \ldots\} \\
c_k &= 2c_{k-1} - c_{k-3} + \frac{\pi^k}{2^k k!} & k &\in \{5, 9, 13, 17, \ldots\} \\
c_k &= 2c_{k-1} - 3c_{k-2} + c_{k-3} & k &\in \{6, 10, 14, 18, \ldots\} \\
c_k &= 3c_{k-2} - 4c_{k-3} - c_{k-4} + c_{k-5} - \frac{\pi^k}{2^k k!} & k &\in \{7, 11, 15, 19, \ldots\}
\end{aligned} \tag{33}$$

By construction the coefficients are positive and decrease according to $0 < c_{k+1} < c_k$, $k \geq 2$, and it is the case that $\varepsilon_1(x_1) \geq 0$ for $x_1 \in [0, 1]$.

**Proof**

The proof is detailed in Appendix 1.

### 4.1.1   Coefficients

The first few coefficients are: $c_2 = 0.14159$, $c_3 = 0.0664249$, $c_4 = 0.057682$, $c_5 = 0.053464$, $c_6 = 3.06823 \times 10^{-4}$, and $c_7 = 1.49925 \times 10^{-4}$.

### 4.1.2   Implication

As $\varepsilon_1(x_1) \geq 0$, $0 \leq x_1 \leq 1$, it follows that $f_1(\pi x_1/2)$ is a lower bound for $\sin(\pi x_1/2)$, $x_1 \in (0, 1)$. The scaling relationships specified in Lemma 3 then imply that $f_1(x)$ is a lower bound for $\sin(x)$, $x \in (0, \pi/2)$.

## 4.2   Lower Bound Based on Second Order Spline Approximation

The second order spline approximation, as defined by Equation 25, is also a lower bound for the sine function. The proof of this is based on considering the scaled error function

$$\varepsilon_2(x_1) = \sin\left[\frac{\pi}{2} \cdot x_1\right] - f_2\left[\frac{\pi}{2} \cdot x_1\right] \qquad x_1 \in [0, 1] \qquad (34)$$

which has the property of being positive on the interval $(0, 1)$ and being such that $\varepsilon_2(0) = \varepsilon_2(1) = 0$:

**Theorem 4.2   Error in Second Order Spline Approximation to Sine**

Using the notation $t = x_1$, the error function defined by Equation 34, can be written as a summation of positive terms according to

$$\begin{aligned}\varepsilon_2(t) &= c_3 t^3(1-t)^3 + c_4 t^4(1-t)^4 + c_5 t^5(1-t)^4 + c_6 t^6(1-t)^3 + c_7 t^7(1-t)^3 + c_8 t^8(1-t)^4 + \ldots \\ &= \sum_{k=3}^{\infty} c_k t^k (1-t)^p \qquad p = \begin{cases} 3 & k = 3, 6, 7, 10, 11, 14, 15, \ldots \\ 4 & k = 4, 5, 8, 9, 12, 13, 16, 17, \ldots \end{cases} \\ &= \sum_{k=3}^{\infty} c_k t^k (1-t)^{\left[\frac{7}{2} + \frac{(-1)^{\lfloor k/2 \rfloor}}{2}\right]}\end{aligned} \qquad (35)$$

The coefficients can be defined iteratively according to

$$c_0 = -1 \qquad c_1 = -2 + \frac{\pi}{2} = 2c_0 + \frac{\pi}{2} \qquad c_2 = -3 + \pi + \frac{\pi^2}{8} = 2c_1 - c_0 + \frac{\pi^2}{2^2 2!} \qquad (36)$$

$$c_3 = -10 + 3\pi + \frac{\pi^2}{8} - \frac{\pi^3}{48} = c_2 + 4c_1 - c_0 - \frac{\pi^3}{2^3 3!} \qquad c_4 = -15 + 5\pi + \frac{\pi^2}{8} - \frac{\pi^3}{16} = 3c_3 - 2c_2 - 4c_1 - c_0 \qquad (37)$$

$$c_5 = -36 + \frac{25\pi}{2} + \frac{\pi^2}{4} - \frac{3\pi^3}{16} + \frac{\pi^5}{2^5 5!} = 3c_4 - c_2 - 3c_1 + \frac{\pi^5}{2^5 5!} \qquad (38)$$

$$c_k = 4c_{k-1} - 6c_{k-2} + c_{k-3} \qquad k \in \{6, 10, 14, 18, \ldots\}$$

$$c_k = 2c_{k-1} - 2c_{k-2} - 2c_{k-3} + c_{k-4} - \frac{\pi^k}{2^k k!} \qquad k \in \{7, 11, 15, 19, \ldots\}$$

$$c_k = 3c_{k-1} - 3c_{k-2} + 4c_{k-3} - c_{k-4} \qquad k \in \{8, 12, 16, 20, \ldots\} \qquad (39)$$

$$c_k = 4c_{k-1} - 3c_{k-2} + c_{k-3} - c_{k-4} + \frac{\pi^k}{2^k k!} \qquad k \in \{9, 13, 17, 21, \ldots\}$$

By construction the coefficients are positive and decrease according to $0 < c_{k+1} < c_k$, $k \geq 3$.

**Proof**

The proof is detailed in Appendix 2.

### 4.2.1    Coefficients

The first few coefficients are: $c_3 = 0.0125144$, $c_4 = 0.00377153$, $c_5 = 0.003325$, $c_6 = 3.18534 \times 10^{-3}$ and $c_7 = 1.02411 \times 10^{-5}$.

### 4.2.2    Implication

As $\varepsilon_2(x_1) \geq 0$, $0 \leq x_1 \leq 1$, it follows that $f_2(\pi x_1/2)$ is a lower bound for $\sin(\pi x_1/2)$, $x_1 \in (0, 1)$. The scaling relationships specified in Lemma 3 then imply that $f_2(x)$ is a lower bound for $\sin(x)$, $x \in (0, \pi/2)$.

## 4.3    Higher Order Spline Approximations as Lower Bounds

The proof that a given higher order spline approximation defined in Theorem 3.1 is also a lower order bound for $\sin(x)$, can be proved in an analogous manner to that detailed in the proofs for Theorem 4.1 and Theorem 4.2.

# 5    Applications

## 5.1    Approximations Consistent with a Set Relative Error Bound

For the case of a set relative error bound of $0.001$ for an approximation to $\sin(x)$, over the interval $[0, \pi/2]$, a second order approximation, as specified by Equation 25, is required. The achieved relative error bound is $3.31 \times 10^{-4}$. The following second order approximation with the approximate coefficients

$$f_2(x) \approx x - \frac{1699 x^3}{10000} + \frac{11 x^4}{2000} + \frac{7 x^5}{1250} \qquad (40)$$

yields a relative error bound of $6.59 \times 10^{-4}$. In comparison, a fifth order Taylor series approximation (Table 2.1) has a relative error bound of $4.52 \times 10^{-3.}$.

To achieve a relative error bound of $10^{-6}$, a fourth order approximation, as specified by Equation 27, is required and the achieved relative error bound is $2.48 \times 10^{-8}$. The following fourth order approximation, with the resolution defined by the given coefficient accuracy, yields a relative error bound of $4.31 \times 10^{-7}$:

$$f_4(x) \approx x - \frac{x^3}{6} + 8.33165 \times 10^{-3} x^5 + 5.17 \times 10^{-6} x^6 - 2.0463 \times 10^{-4} x^7 +$$
$$3.55 \times 10^{-6} x^8 + 1.89 \times 10^{-6} x^9 \qquad (41)$$

A ninth order Taylor series approximation (Table 2.1) yields a relative error bound of $3.54 \times 10^{-6}$.

For relative error bounds of $10^{-10}$ and $10^{-16}$, spline approximations of orders six and eight are required. The required approximations are:

$$f_6(x) \approx x - \frac{x^3}{6} + \frac{x^5}{120} - 1.98412876 \times 10^{-4} x^7 + 7.78 \times 10^{-10} x^8 + 2.754282 \times 10^{-6} x^9 +$$
$$1.485 \times 10^{-9} x^{10} - 2.5944 \times 10^{-8} x^{11} + 3.06 \times 10^{-10} x^{12} + 1.11 \times 10^{-10} x^{13} \quad (42)$$

$$f_8(x) \approx x - \frac{x^3}{6} + \frac{x^5}{120} - \frac{x^7}{5040} + 2.755731916334 \times 10^{-6} x^9 +$$
$$3.43891 \times 10^{-14} x^{10} - 2.50521948 \times 10^{-8} x^{11} + 1.26122 \times 10^{-13} x^{12} + 1.604729071 \times 10^{-10} x^{13} + \quad (43)$$
$$7.21759 \times 10^{-14} x^{14} - 7.936185 \times 10^{-13} x^{15} + 7.0722 \times 10^{-15} x^{16} + 1.956 \times 10^{-15} x^{17}$$

The achieved relative error bounds, respectively, are $4.39 \times 10^{-11}$ and $5.66 \times 10^{-17}$. Higher precision in the coefficients leads to relative error bounds, respectively, of $3.91 \times 10^{-13}$ and $2.02 \times 10^{-18}$. For comparison, a thirteenth and a seventeenth order Taylor series approximation have, respectively, relative error bounds of $6.63 \times 10^{-10}$ and $4.35 \times 10^{-14}$.

## 5.2    Lower Bounds for Sine

**Theorem 5.1   Lower Bounds for Sine**

The spline based series for $\sin(x)$, defined in Theorem 3.1, are lower bounds for $\sin(x)$ on the interval $[0, \pi/2]$. Thus, for example, the first and second order spline approximations detailed in Theorem 3.1, yield the bounds:

$$0 < 1 + \frac{12}{\pi^2}\left[1 - \frac{\pi}{3}\right]x - \frac{16}{\pi^3}\left[1 - \frac{\pi}{4}\right]x^2 < \frac{\sin(x)}{x} \qquad x \in (0, \pi/2) \quad (44)$$

$$0 < 1 + \frac{80}{\pi^3}\left[1 - \frac{3\pi}{10} - \frac{\pi^2}{80}\right]x^2 - \frac{240}{\pi^4}\left[1 - \frac{4\pi}{15} - \frac{\pi^2}{60}\right]x^3 + \frac{192}{\pi^5}\left[1 - \frac{\pi}{4} - \frac{\pi^2}{48}\right]x^4 < \frac{\sin(x)}{x} \qquad x \in (0, \pi/2) \quad (45)$$

and similar results hold for higher order spline approximations. There is equality at the points $0, \pi/2$. The maximum relative errors associated with these lower bounds are detailed in Table 3.1.

**Proof**

This result follows from Theorem 4.1 and Theorem 4.2 and, similarly, for higher order spline approximations.

## 5.3    New Series for Sine: First Order

Consider the result stated in Theorem 4.1 which specifies a series for the error function $\varepsilon_1$. Convergence of this series implies a new series for the sine function:

**Theorem 5.2   Series Approximation for Sin Based on First Order Spline Approximation**

The error function associated with a first order spline approximation yields the following series for the sine function that is valid for the interval $[0, \pi/2]$:

$$\sin(x) = \frac{2x}{\pi} + \frac{2x}{\pi}\left[1 - \frac{2x}{\pi}\right] + \sum_{k=1}^{\infty} c_k \left[\frac{2x}{\pi}\right]^k \left[1 - \frac{2x}{\pi}\right]^p \qquad p = \begin{cases} 2 & k = 1, 2, 5, 6, 9, 10, 13, 14, 17, \ldots \\ 3 & k = 3, 4, 7, 8, 11, 12, 15, 16, 19, \ldots \end{cases} \quad (46)$$

where $c_1 = 2[-1 + \pi/4]$ and the coefficients $c_2, c_3, \ldots$ are defined in Theorem 4.1.

The first few terms are

$$\sin(x) = \frac{2x}{\pi} + \frac{2x}{\pi}\left[1 - \frac{2x}{\pi}\right] + 2\left[-1 + \frac{\pi}{4}\right]\frac{2x}{\pi}\left[1 - \frac{2x}{\pi}\right]^2 + (\pi - 3)\left[\frac{2x}{\pi}\right]^2\left[1 - \frac{2x}{\pi}\right]^2 +$$
$$\left[-4 + \frac{3\pi}{2} - \frac{\pi^3}{48}\right]\left[\frac{2x}{\pi}\right]^3\left[1 - \frac{2x}{\pi}\right]^3 + \left[-9 + \frac{7\pi}{2} - \frac{\pi^3}{16}\right]\left[\frac{2x}{\pi}\right]^4\left[1 - \frac{2x}{\pi}\right]^3 + \ldots$$

(47)

**Proof**

The proof is detailed in Appendix 3.

### 5.3.1   Approximations and Relative Errors

The sine function can be approximated by taking the first $n$ terms in the series defined in Theorem 5.2, i.e.

$$\sin(x) \approx \frac{2x}{\pi} + \frac{2x}{\pi}\left[1 - \frac{2x}{\pi}\right] + \sum_{k=1}^{n} c_k\left[\frac{2x}{\pi}\right]^k\left[1 - \frac{2x}{\pi}\right]^p \qquad p = \begin{cases} 2 & k = 1, 2, 5, 6, 9, 10, \ldots \\ 3 & k = 3, 4, 7, 8, 11, 12, \ldots \end{cases}$$

(48)

The graphs of the relative errors associated with the approximations to $\sin(x)$, which are defined by the first to ninth terms in the series approximations, are shown in Figure 5. The associated relative error bounds are detailed in Table 5.1.

Table 5.1   Relative error bounds for the first and second order series approximations to $\sin(x)$.

| # terms: $n$ | Relative error bound. First order: Equation 48. | Relative error bound. Second order: Equation 52. |
|---|---|---|
| 1  | $1.63 \times 10^{-2}$  |  |
| 2  | $2.70 \times 10^{-3}$  | $3.31 \times 10^{-4}$  |
| 3  | $1.42 \times 10^{-3}$  | $3.89 \times 10^{-5}$  |
| 4  | $9.14 \times 10^{-4}$  | $2.00 \times 10^{-5}$  |
| 6  | $1.30 \times 10^{-6}$  | $3.11 \times 10^{-8}$  |
| 8  | $7.92 \times 10^{-7}$  | $3.91 \times 10^{-9}$  |
| 12 | $1.50 \times 10^{-10}$ | $2.83 \times 10^{-13}$ |
| 16 | $9.85 \times 10^{-15}$ | $9.05 \times 10^{-18}$ |
| 20 | $2.82 \times 10^{-19}$ | $1.45 \times 10^{-22}$ |

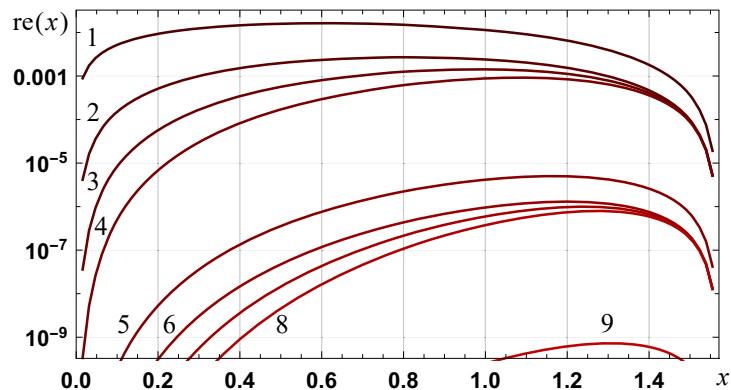

Figure 5   Graphs of the relative errors in the approximation to $\sin(x)$, as defined by Equation 48, for the cases of the first to ninth terms.

## 5.4    New Series for Sine: Second Order

The series for the error function $\varepsilon_2$ defined in Theorem 4.2 converges and this underpins the definition of the following series for the sine function:

### Theorem 5.3   Series Approximation for Sin Based on Second Order Spline Approximation

The error function associated with a second order spline approximation yields the following series for the interval $[0, \pi/2]$:

$$\sin(x) = 1 - \frac{\pi^2}{8}\left[1 - \frac{2x}{\pi}\right]^2 + \sum_{k=0}^{\infty} c_k \left[\frac{2x}{\pi}\right]^k \left[1 - \frac{2x}{\pi}\right]^p \qquad p = \begin{cases} 3 & k = 2, 3, 6, 7, 10, 11, 14, 15, \ldots \\ 4 & k = 0, 1, 4, 5, 8, 9, 12, 13, 16, 17, \ldots \end{cases} \qquad (49)$$

where

$$c_0 = -1 + \frac{\pi^2}{8} \qquad c_1 = -4 + \frac{\pi}{2} + \frac{\pi^2}{4} \qquad c_2 = -10 + 2\pi + \frac{3\pi^2}{8} \qquad (50)$$

and with the coefficients $c_3, c_4, \ldots$ being defined in Theorem 4.2. The first few terms are

$$\begin{aligned}\sin(x) = {} & 1 - \frac{\pi^2}{8}\left[1 - \frac{2x}{\pi}\right]^2 + \left[-1 + \frac{\pi^2}{8}\right]\left[1 - \frac{2x}{\pi}\right]^4 + \left[-4 + \frac{\pi}{2} + \frac{\pi^2}{4}\right]\left[\frac{2x}{\pi}\right]\left[1 - \frac{2x}{\pi}\right]^4 + \\ & \left[-10 + 2\pi + \frac{3\pi^2}{8}\right]\left[\frac{2x}{\pi}\right]^2\left[1 - \frac{2x}{\pi}\right]^3 + \left[-10 + 3\pi + \frac{\pi^2}{8} - \frac{\pi^3}{48}\right]\left[\frac{2x}{\pi}\right]^3\left[1 - \frac{2x}{\pi}\right]^3 + \ldots\end{aligned} \qquad (51)$$

**Proof**

The proof is detailed in Appendix 4.

### 5.4.1    Approximations and Relative Errors

The sine function can be approximated by taking the first $n$ terms in the series defined in Theorem 5.3, i.e.

$$\sin(x) \approx 1 - \frac{\pi^2}{8}\left[1 - \frac{2x}{\pi}\right]^2 + \sum_{k=0}^{n} c_k \left[\frac{2x}{\pi}\right]^k \left[1 - \frac{2x}{\pi}\right]^p \qquad p = \begin{cases} 3 & k = 2, 3, 6, 7, 10, 11, \ldots \\ 4 & k = 0, 1, 4, 5, 8, 9, \ldots \end{cases} \qquad (52)$$

The graphs of the relative errors in the approximations, for the case of two to nine terms, are shown in Figure 6. The associated relative error bounds are detailed in Table 5.1.

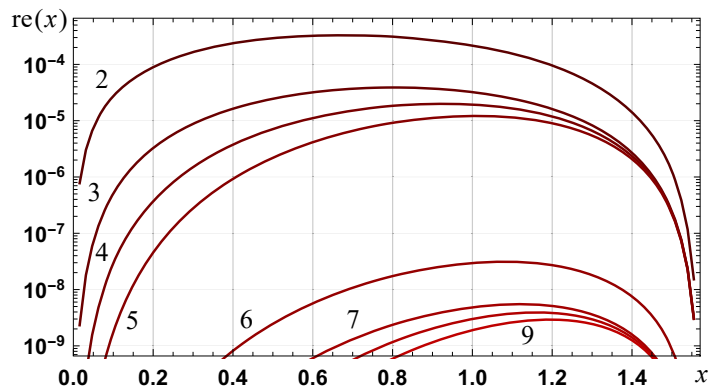

**Figure 6**   Graphs of the relative errors in the approximations to $\sin(x)$ as defined by Equation 52, for orders two to nine.

### 5.4.2    New Series for Sine: Higher Order Series

Higher order series for $\sin(x)$ can be generated in a similar manner and, as the results in Table 5.1 indicate, with lower relative error bounds.

## 5.5    Comparison of Series for Sine

The approximations for $\sin(x)$, specified by Theorem 3.1, Theorem 5.2 and Theorem 5.3, of order two, are:

$$f_2(x) = x + \frac{80}{\pi^3}\left[1 - \frac{3\pi}{10} - \frac{\pi^2}{80}\right]x^3 - \frac{240}{\pi^4}\left[1 - \frac{4\pi}{15} - \frac{\pi^2}{60}\right]x^4 + \frac{192}{\pi^5}\left[1 - \frac{\pi}{4} - \frac{\pi^2}{48}\right]x^5 \quad (53)$$

$$f_2(x) = \frac{2x}{\pi} + \frac{2x}{\pi}\left[1 - \frac{2x}{\pi}\right] + 2\left[-1 + \frac{\pi}{4}\right]\frac{2x}{\pi}\left[1 - \frac{2x}{\pi}\right]^2 + (\pi - 3)\left[\frac{2x}{\pi}\right]^2\left[1 - \frac{2x}{\pi}\right]^2 \quad (54)$$

$$f_2(x) = 1 - \frac{\pi^2}{8}\left[1 - \frac{2x}{\pi}\right]^2 + \left[-1 + \frac{\pi^2}{8}\right]\left[1 - \frac{2x}{\pi}\right]^4 + \left[-4 + \frac{\pi}{2} + \frac{\pi^2}{4}\right]\left[\frac{2x}{\pi}\right]\left[1 - \frac{2x}{\pi}\right]^4 +$$

$$\left[-10 + 2\pi + \frac{3\pi^2}{8}\right]\left[\frac{2x}{\pi}\right]^2\left[1 - \frac{2x}{\pi}\right]^3 \quad (55)$$

with respective relative error bounds over the interval $[0, \pi/2]$ of $3.31 \times 10^{-4}$, $2.70 \times 10^{-3}$ and $3.31 \times 10^{-4}$. In general, the series specified by Theorem 3.1 have better convergence that the other two series and a simpler form.

## 5.6    Upper Bounds for Sine

Consistent with Lemma 1, the requirement for an upper bound for the sine function, based on the first and second order spline approximations $f_1$ and $f_2$, as defined by Equation 24 and Equation 25, is for

$$\varepsilon_1\left[\frac{2x}{\pi}\right] > 2\varepsilon_2\left[\frac{2x}{\pi}\right] \qquad x \in \left(0, \frac{\pi}{2}\right) \quad (56)$$

where $\varepsilon_1$ and $\varepsilon_2$ are defined by Equation 30 and Equation 35. When this is the case the following upper bound can be defined:

**Theorem 5.4    Upper Bound for Sine**

The lower bounds $f_1$ and $f_2$, as defined by Equation 24 and Equation 25, define an upper bound $f_2^U$ for $\sin(x)$ according to

$$f_2^U(x) = 2f_2(x) - f_1(x) \qquad x \in (0, \pi/2) \quad (57)$$

i.e.

$$1 - \frac{12}{\pi^2}\left[1 - \frac{\pi}{3}\right]x + \frac{176}{\pi^3}\left[1 - \frac{13\pi}{44} - \frac{\pi^2}{88}\right]x^2 - \frac{480}{\pi^4}\left[1 - \frac{4\pi}{15} - \frac{\pi^2}{60}\right]x^3 + \frac{384}{\pi^5}\left[1 - \frac{\pi}{4} - \frac{\pi^2}{48}\right]x^4 > \frac{\sin(x)}{x}$$

$$x \in (0, \pi/2) \quad (58)$$

The error in this bound, as defined by

$$\varepsilon_2^U(x) = f_2^U(x) - \sin(x) \qquad x \in (0, \pi/2) \quad (59)$$

is shown in Figure 7. Note that the error $\varepsilon_2^U(x)$ is close to the error in the lower bound, $\varepsilon_1(x)$, associated with $f_1(x)$, which is shown in Figure 4.

**Proof**
The proof is detailed in Appendix 5.

### 5.6.1    Higher Order Approximations

Higher order approximations for the upper bound of the sine function follow in an analogous manner. For example

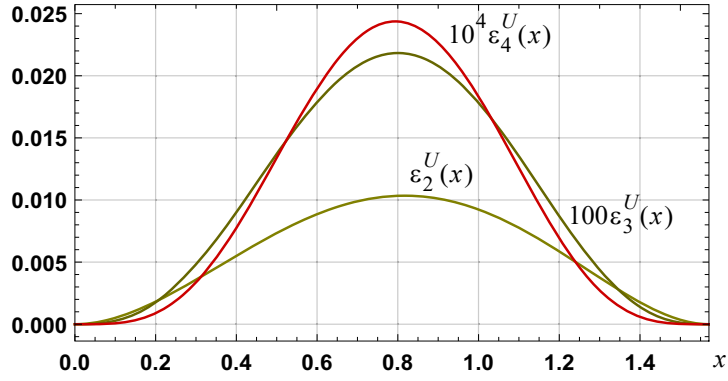

**Figure 7**   Graphs of the errors in the functions $f_2^U, f_3^U$ and $f_4^U$, which defined second, third and fourth order upper bounds to $\sin(x)$.

$$f_3^U(x) = 2f_3(x) - f_2(x) \qquad x \in (0, \pi/2)$$
$$f_4^U(x) = 2f_4(x) - f_3(x) \qquad x \in (0, \pi/2)$$
(60)

are valid upper bounds with $f_2, f_3$ and $f_4$ being defined by Equation 25 to Equation 27. The graphs of the errors in these bounds, denoted, respectively, $\varepsilon_3^U$ and $\varepsilon_4^U$, are shown in Figure 7.

### 5.7   Bounds for Cosine

The results detailed in Lemma 2 allow the bounds for the cosine function to be defined by utilizing the bounds specified for the sine function:

**Theorem 5.5   Lower Bounds for Cosine**

First to fourth order spline based approximations for $\sin(x)$ over the interval $[0, \pi/2]$, as specified in Theorem 3.1, lead to the following lower bounds for the cosine function:

$$g_1(y) = 1 - \frac{12}{\pi^2}\left[1 - \frac{\pi}{6}\right]y^2 + \frac{16}{\pi^3}\left[1 - \frac{\pi}{4}\right]y^3 \tag{61}$$

$$g_2(y) = 1 - \frac{y^2}{2} - \frac{80}{\pi^3}\left[1 - \frac{\pi}{5} - \frac{3\pi^2}{80}\right]y^3 + \frac{240}{\pi^4}\left[1 - \frac{7\pi}{30} - \frac{\pi^2}{40}\right]y^4 - \frac{192}{\pi^5}\left[1 - \frac{\pi}{4} - \frac{\pi^2}{48}\right]y^5 \tag{62}$$

$$g_3(y) = 1 - \frac{y^2}{2} - \frac{560}{\pi^4}\left[1 - \frac{3\pi}{14} - \frac{\pi^2}{28} + \frac{\pi^3}{1680}\right]y^4 + \frac{2688}{\pi^5}\left[1 - \frac{13\pi}{56} - \frac{5\pi^2}{168} + \frac{\pi^3}{1344}\right]y^5 -$$
$$\frac{4480}{\pi^6}\left[1 - \frac{17\pi}{70} - \frac{3\pi^2}{112} + \frac{\pi^3}{1120}\right]y^6 + \frac{2560}{\pi^7}\left[1 - \frac{\pi}{4} - \frac{\pi^2}{40} + \frac{\pi^3}{960}\right]y^7 \tag{63}$$

$$g_4(y) = 1 - \frac{y^2}{2} + \frac{y^4}{24} - \frac{4032}{\pi^5}\left[1 - \frac{2\pi}{9} - \frac{5\pi^2}{144} + \frac{\pi^3}{1008} + \frac{5\pi^4}{48384}\right]y^5 +$$
$$\frac{26880}{\pi^6}\left[1 - \frac{7\pi}{30} - \frac{\pi^2}{32} + \frac{23\pi^3}{20160} + \frac{\pi^4}{16128}\right]y^6 - \frac{69120}{\pi^7}\left[1 - \frac{13\pi}{54} - \frac{7\pi^2}{240} + \frac{11\pi^3}{8640} + \frac{\pi^4}{20736}\right]y^7 +$$
$$\frac{80640}{\pi^8}\left[1 - \frac{31\pi}{126} - \frac{\pi^2}{36} + \frac{\pi^3}{720} + \frac{\pi^4}{24192}\right]y^8 - \frac{35840}{\pi^9}\left[1 - \frac{\pi}{4} - \frac{3\pi^2}{112} + \frac{\pi^3}{672} + \frac{\pi^4}{26880}\right]y^9 \tag{64}$$

**Proof**

These results follow from the result stated in Lemma 2 of $f_L(\pi/2 - y) < \cos(y)$, $y \in [0, \pi/2]$ and the spline approximations for the sine function stated in Theorem 3.1. Together these imply $g_k(y) < \cos(y)$ where $g_k(y) = f_k(\pi/2 - y)$, $k \in \{1, 2, 3, 4\}$. Alternatively, they can be directly derived from Equation 21 for the case of $f(x) = \cos(x)$.

### 5.7.1   Results

The relative error bounds in the approximations for the cosine function defined in Theorem 5.5 are identical to the relative error bounds associated with the lower bounds for the sine function which are detailed in Table 3.1.

### 5.7.2   Upper Bounds

Upper bounds for the cosine function can be generated consistent with the results stated in Theorem 5.4. For example, the result $g_2^U(y) = 2g_2(y) - g_1(y)$, $y \in (0, \pi/2)$ implies

$$1 + \frac{12}{\pi^2}\left[1 - \frac{\pi}{6} - \frac{\pi^2}{12}\right]y^2 - \frac{176}{\pi^3}\left[1 - \frac{9\pi}{44} - \frac{3\pi^2}{88}\right]y^3 + \frac{480}{\pi^4}\left[1 - \frac{7\pi}{30} - \frac{\pi^2}{40}\right]y^4 - \frac{384}{\pi^5}\left[1 - \frac{\pi}{4} - \frac{\pi^2}{48}\right]y^5 \tag{65}$$

$$> \cos(y) \qquad y \in (0, \pi/2)$$

The error in this bound is

$$\varepsilon_{2,\cos}^U(y) = g_2^U(y) - \cos(y) = 2g_2(y) - g_1(y) - \cos(y)$$

$$= 2f_2\left[\frac{\pi}{2} - y\right] - f_1\left[\frac{\pi}{2} - y\right] - \sin\left[\frac{\pi}{2} - y\right] = \varepsilon_2^U\left[\frac{\pi}{2} - y\right] \tag{66}$$

and, thus, has the same bound as $\varepsilon_2^U$ whose graph is shown in Figure 7. Higher order upper bounds can be defined in a similar manner.

### 5.7.3   Comparison with Published Bounds

The list of inequalities for $\cos(y)$ over the interval $[0, \pi/2]$ is similar in scope to that for $\sin(x)$ with the Kober inequality (Kober 1944) serving as a reference inequality. Qi et. al. (2009) provides a useful overview. In general, the bounds proposed in this paper are closer bounds than published results and can be made arbitrary accurate.

## 5.8   Bounds for the Sine Integral

Bounds for $\sin(x)/x$ readily lead to bounds for the sine integral function (e.g. Lv et al. 2017, Proposition 5; Zeng & Wu 2013, Theorem 9) which is defined according to

$$\text{Si}(x) = \int_0^x \frac{\sin(\lambda)}{\lambda} d\lambda \tag{67}$$

The lower bound published by Lv 2017 is

$$\frac{2x + \sin(x)}{3} - \frac{x^3 + 3x\cos(x) - 3\sin(x)}{9\pi^2} < \text{Si}(x) \tag{68}$$

Integration of the spline series for $\sin(x)$ defined by Theorem 3.1 (after dividing by $x$) leads to the following approximations:

**Theorem 5.6   Spline Based Approximations for Sine Integral**

First to fourth order approximations for the sine integral function, which are also lower bound functions, are:

$$h_1(x) = x + \frac{6}{\pi^2}\left[1 - \frac{\pi}{3}\right]x^2 - \frac{16}{3\pi^3}\left[1 - \frac{\pi}{4}\right]x^3 \tag{69}$$

$$h_2(x) = x + \frac{80}{3\pi^2}\left[1 - \frac{3\pi}{10} - \frac{\pi^2}{80}\right]x^3 - \frac{60}{\pi^4}\left[1 - \frac{4\pi}{15} - \frac{\pi^2}{60}\right]x^4 + \frac{192}{5\pi^5}\left[1 - \frac{\pi}{4} - \frac{\pi^2}{48}\right]x^5 \tag{70}$$

$$h_3(x) = x - \frac{x^3}{18} + \frac{140}{\pi^4}\left[1 - \frac{2\pi}{7} - \frac{\pi^2}{56} + \frac{\pi^3}{420}\right]x^4 - \frac{2688}{5\pi^5}\left[1 - \frac{15\pi}{56} - \frac{\pi^2}{48} + \frac{\pi^3}{672}\right]x^5 +$$

$$\frac{2240}{3\pi^6}\left[1 - \frac{9\pi}{35} - \frac{13\pi^2}{560} + \frac{\pi^3}{840}\right]x^6 - \frac{2560}{7\pi^7}\left[1 - \frac{\pi}{4} - \frac{\pi^2}{40} + \frac{\pi^3}{960}\right]x^7 \tag{71}$$

$$h_4(x) = x - \frac{x^3}{18} + \frac{4032}{5\pi^5}\left[1 - \frac{5\pi}{18} - \frac{\pi^2}{48} + \frac{5\pi^3}{2016} + \frac{\pi^4}{48384}\right]x^5 -$$

$$\frac{4480}{\pi^6}\left[1 - \frac{4\pi}{15} - \frac{11\pi^2}{480} + \frac{\pi^3}{504} + \frac{\pi^4}{40320}\right]x^6 + \frac{69120}{7\pi^7}\left[1 - \frac{7\pi}{27} - \frac{53\pi^2}{2160} + \frac{\pi^3}{576} + \frac{\pi^4}{34560}\right]x^7 - \tag{72}$$

$$\frac{10080}{\pi^8}\left[1 - \frac{16\pi}{63} - \frac{13\pi^2}{504} + \frac{\pi^3}{630} + \frac{\pi^4}{30240}\right]x^8 + \frac{35840}{9\pi^9}\left[1 - \frac{\pi}{4} - \frac{3\pi^2}{112} + \frac{\pi^3}{672} + \frac{\pi^4}{26880}\right]x^9$$

Higher order approximations can similarly be generated.

### Proof

First to fourth order spline based approximations for the sine integral, and for the interval $[0, \pi/2]$, can be generated from the approximations detailed in Theorem 3.1 for the sine function according to

$$h_k(x) = \int_0^x \frac{f_k(\lambda)}{\lambda} d\lambda \qquad k \in \{0, 1, \ldots\} \tag{73}$$

The stated results then follow from the definitions for $f_k$ given in Theorem 3.1. As $f_k$, $k \in \{0, 1, \ldots\}$ are lower bounds for the sine function, it then follows that the approximations for the sine integral function are also lower bounds.

### 5.8.1    Results

Graphs of the relative error in the approximations defined in Theorem 5.6 are shown in Figure 8 and show the improvement over the lower bound proposed by Lv 2017 for orders three and higher. Bounds on the relative error are tabulated in Table 5.2 and clearly can be made arbitrarily small by increasing the order of approximation.

Table 5.2    Relative error bounds for the spline based approximations to the sine integral function over the interval $[0, \pi/2]$.

| Order of approx. | Relative error bound |
| --- | --- |
| 0 | 0.363 |
| 1 | $1.24 \times 10^{-2}$ |
| 2 | $2.12 \times 10^{-4}$ |
| 3 | $2.06 \times 10^{-6}$ |
| 4 | $1.28 \times 10^{-8}$ |
| 8 | $8.21 \times 10^{-19}$ |
| 12 | $1.61 \times 10^{-30}$ |
| 16 | $2.85 \times 10^{-43}$ |

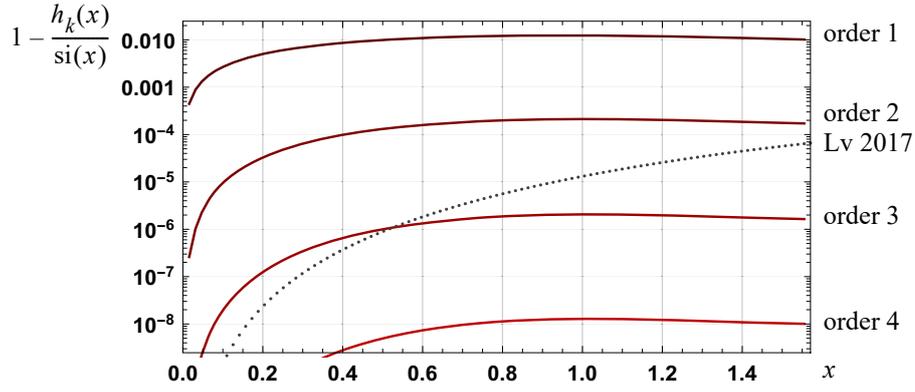

**Figure 8** Graph of the relative error in spline based approximations, of orders one to four, for the sine integral function as well as the relative error in the lower bound specified by Lv (2017) and stated in Equation 68.

### 5.8.2 Notes

Upper bound functions for the sine integral can be generated from the lower bounds and consistent with Lemma 1.

## 6 Conclusion

Spline based approximations for $\sin(x)$ over the interval $[0, \pi/2]$ have been to detailed and exhibit better convergence that comparative Taylor series approximations. It was proved that these approximations are also lower bounds for $\sin(x)$ and $\sin(x)/x$, over the interval $[0, \pi/2]$, with increasingly high accuracy and with lower relative error bounds than published lower bounded functions. For example, approximations of order two, four and eight have relative error bounds, respectively, of $3.31 \times 10^{-4}$, $2.48 \times 10^{-8}$ and $2.02 \times 10^{-18}$.

The analysis underpinning the proof that the approximations are lower bounds for the sine function led to new series for $\sin(x)$. These have better convergence that Taylor series for $\sin(x)$ over the interval $[0, \pi/2]$.

Applications were detailed and include, first, a sequence of upper bounded functions, with increasing accuracy, for $\sin(x)$ over the interval $[0, \pi/2]$. Second, a sequence of upper and lower bounded functions for $\cos(x)$ with accuracy consistent with the corresponding sine approximations. Finally, a sequence of convergent lower bounded approximations for the sine integral function over the interval $[0, \pi/2]$.

**Acknowledgement** The support of Prof. A. Zoubir, SPG, Technische Universität Darmstadt, Darmstadt, Germany, who hosted a visit where the writing of this paper was completed, is gratefully acknowledged.

## Appendix 1: Proof of Theorem 4.1

By utilizing a Taylor series expansion for $\sin(\pi t/2)$, $t \in [0, 1]$, and the definition for $f_1$ defined by Equation 24, the error in a first order spline function approximation for $\sin(\pi t/2)$ is

$$\varepsilon_1(t) = \sin\left[\frac{\pi}{2} \cdot t\right] - f_1\left[\frac{\pi}{2} \cdot t\right]$$

$$= (\pi - 3)t^2 + \left[2 - \frac{\pi}{2} - \frac{\pi^3}{48}\right]t^3 + \frac{1}{5!} \cdot \left[\frac{\pi}{2} \cdot t\right]^5 - \frac{1}{7!} \cdot \left[\frac{\pi}{2} \cdot t\right]^7 + \frac{1}{9!} \cdot \left[\frac{\pi}{2} \cdot t\right]^9 + \ldots \quad (74)$$

$$= e_2 t^2 + e_3 t^3 + e_5 t^5 + e_7 t^7 + e_9 t^9 + \ldots$$

where $e_2 = \pi - 3 = 0.314159$ and $e_3 = 2 - \pi/2 - \pi^3/48 = -0.21676$ etc. The goal is to write $\varepsilon_1$ as a summation of positive terms only. To this end the following form is considered:

$$\varepsilon_1(t) = c_2 t^2 (1-t)^{p_2} + c_3 t^3 (1-t)^{p_3} + c_4 t^4 (1-t)^{p_4} + \ldots + c_n t^n (1-t)^{p_n} + \ldots$$

$$= c_2 t^2 \left[1 - p_2 t + \binom{p_2}{2} t^2 - \ldots \pm t^{p_2}\right] + c_3 t^3 [1 - p_3 t + \ldots \pm t^{p_3}] + c_4 t^4 [1 - p_4 t + \ldots \pm t^{p_4}] + \ldots \quad (75)$$

$$= c_2 t^2 + [-p_2 c_2 + c_3] t^3 + \left[\binom{p_2}{2} c_2 - p_3 c_3 + c_4\right] t^4 + \ldots$$

Equating Equation 74 and Equation 75 it immediately follows that $c_2 = e_2 > 0$. For $c_3$ to be positive, the factor $p_2$ needs to be chosen such that $e_3 = -p_2 c_2 + c_3$ and $c_3 = e_3 + p_2 c_2 > 0$. This is guaranteed when $p_2 > -e_3/c_2$. Choosing the smallest integer greater than zero such that this is the case yields $p_2 = \lfloor -e_3/e_2 \rfloor + 1$ and for the given values of $e_2$ and $e_3$ the required value is $p_2 = 2$. It then follows that $c_3 = -4 + 3\pi/2 - \pi^3/48 = 0.0664249$.

The ratio of $c_3$ to $c_2$ then is

$$\frac{c_3}{c_2} = \frac{e_3 + p_2 c_2}{e_2} = \frac{e_3}{e_2} + p_2 = \frac{e_3}{e_2} + \left\lfloor \frac{-e_3}{e_2} \right\rfloor + 1 = 1 - \left[\frac{-e_3}{e_2} - \left\lfloor \frac{-e_3}{e_2} \right\rfloor\right] \quad (76)$$

which implies $0 < c_3/c_2 < 1$ as $-e_3/e_2 > 0$ and $0 \le x - \lfloor x \rfloor < 1$.

With $p_2 = 2$ and with $c_2$ and $c_3$ set, it follows that

$$\varepsilon_1(t) = c_2 t^2 [1 - 2t + t^2] + c_3 t^3 \left[1 - p_3 t + \binom{p_3}{2} t^2 - \ldots \pm t^{p_3}\right] + c_4 t^4 [1 - p_4 t + \ldots \pm t^{p_4}]$$

$$= e_2 t^2 + e_3 t^3 + [c_2 - p_3 c_3 + c_4] t^4 + c_3 t^3 \left[\binom{p_3}{2} t^2 - \ldots \pm t^{p_3}\right] + c_4 t^4 [-p_4 t^2 + \ldots \pm t^{p_4}] + \ldots \quad (77)$$

and the requirement of $e_4 = c_2 - p_3 c_3 + c_4 = 0$ can be satisfied by choosing $p_3$ such that $c_4 > 0$. With $c_4 = -c_2 + p_3 c_3 > 0$ the requirement is for $p_3 > c_2/c_3$. Choosing the smallest integer greater than zero such that this is the case yields $p_3 = \lfloor c_2/c_3 \rfloor + 1$ and for the defined values of $c_2$ and $c_3$ it is the case that $p_3 = 3$. It then follows that $c_4 = -9 + 7\pi/2 - \pi^3/16 = 0.057682$.

The ratio of $c_4$ to $c_3$ then is

$$\frac{c_4}{c_3} = \frac{-c_2 + p_3 c_3}{c_3} = \frac{-c_2}{c_3} + p_3 = 1 + \left[\left\lfloor \frac{c_2}{c_3} \right\rfloor - \frac{c_2}{c_3}\right] \quad (78)$$

which implies $0 < c_4/c_3 < 1$.

Proceeding in this manner yields the stated form and proof of positive coefficients with decreasing values. The coefficients are defined according to

$$c_2 = -3 + \pi \qquad c_3 = -4 + \frac{3\pi}{2} - \frac{\pi^3}{48} \quad (79)$$

$$c_4 = -9 + \frac{7\pi}{2} - \frac{\pi^3}{16} \qquad c_5 = -15 + 6\pi - \frac{\pi^3}{8} + \frac{\pi^5}{3840} \quad (80)$$

$$c_6 = -7 + 3\pi - \frac{\pi^3}{12} + \frac{\pi^5}{1920} \qquad c_7 = -8 + \frac{7\pi}{2} - \frac{5\pi^3}{48} + \frac{\pi^5}{1280} - \frac{\pi^7}{645120} \quad (81)$$

$$c_8 = -17 + \frac{15\pi}{2} - \frac{11\pi^3}{48} + \frac{7\pi^5}{3840} - \frac{\pi^7}{215040} \qquad c_9 = -27 + 12\pi - \frac{3\pi^3}{8} + \frac{\pi^5}{320} - \frac{\pi^7}{107520} + \frac{\pi^9}{185,794,560} \quad (82)$$

$$c_{10} = -11 + 5\pi - \frac{\pi^3}{6} + \frac{\pi^5}{640} - \frac{\pi^7}{161280} + \frac{\pi^9}{92,897,280} \tag{83}$$

$$c_{11} = -12 + \frac{11\pi}{2} - \frac{3\pi^3}{16} + \frac{7\pi^5}{3840} - \frac{\pi^7}{129024} + \frac{\pi^9}{61,931,520} - \frac{\pi^{11}}{81,749,606,400} \tag{84}$$

$$c_{12} = -25 + \frac{23\pi}{2} - \frac{19\pi^3}{48} + \frac{\pi^5}{256} - \frac{11\pi^7}{645120} + \frac{\pi^9}{26,542,080} - \frac{\pi^{11}}{27,249,868,800} \tag{85}$$

To find a general formula for the coefficients consider the definitions:

$$c_0 = -1 \qquad c_1 = -2 + \frac{\pi}{2} = 2c_0 + \frac{\pi}{2} \tag{86}$$

It then can be shown that

$$c_2 = 2c_1 - c_0 \qquad c_3 = c_2 + c_1 - c_0 - \frac{\pi^3}{2^3 3!} \tag{87}$$

$$c_4 = 3c_3 - c_2 \qquad c_5 = 2c_4 - c_2 + \frac{\pi^5}{2^5 5!} \tag{88}$$

$$c_6 = 2c_5 - 3c_4 + c_3 \qquad c_7 = 3c_5 - 4c_4 - c_3 + c_2 - \frac{\pi^7}{2^7 7!} \tag{89}$$

The iteration formulas

$$\begin{aligned} c_k &= 3c_{k-1} - c_{k-2} & k \in \{4, 8, 12, 16, \ldots\} \\ c_k &= 2c_{k-1} - c_{k-3} + \frac{\pi^k}{2^k k!} & k \in \{5, 9, 13, 17, \ldots\} \\ c_k &= 2c_{k-1} - 3c_{k-2} + c_{k-3} & k \in \{6, 10, 14, 18, \ldots\} \\ c_k &= 3c_{k-2} - 4c_{k-3} - c_{k-4} + c_{k-5} - \frac{\pi^k}{2^k k!} & k \in \{7, 11, 15, 19, \ldots\} \end{aligned} \tag{90}$$

then follow.

## A1.1  Proof of Convergence

To prove that the new series for $\varepsilon_1(t)$, as specified by Equation 75, converges over the interval $t \in [0, 1]$, consider the upper bounds:

$$\begin{aligned} \varepsilon_1(t) &= c_2 t^2 (1-t)^2 + c_3 t^3 (1-t)^3 + c_4 t^4 (1-t)^3 + c_5 t^5 (1-t)^2 + c_6 t^6 (1-t)^2 + c_7 t^7 (1-t)^3 + \ldots \\ &< c_2 t^2 + c_3 t^3 + c_4 t^4 + c_5 t^5 + c_6 t^6 + c_7 t^7 + \ldots \\ &< c_2 t^2 [1 + t + t^2 + t^3 + t^4 + t^5 + \ldots] = c_2 t^2 \cdot \frac{1}{1-t} \qquad t \in (0, 1) \end{aligned} \tag{91}$$

The inequalities follow as $c_k > 0$, $0 < c_{k+1}/c_k < 1$ for all $k$ and the final form, consistent with the convergence of a geometric series, is valid for $t \in (0, 1)$. As $\varepsilon_1(0) = \varepsilon_1(1) = 0$ it follows that the series defining $\varepsilon_1(t)$ is bounded above on the interval $t \in [0, 1]$.

With positive coefficients, the series defining $\varepsilon_1(t)$ is bounded below according to $\varepsilon_1(t) \geq 0$, $0 \leq t \leq 1$.

## Appendix 2:    Proof of Theorem 4.2

By utilizing a Taylor series expansion for $\sin[\pi t/2]$, $t \in [0, 1]$, and the definition for $f_2$ defined by Equation 25, the error in a second order spline function approximation for $\sin[\pi t/2]$ is

$$\varepsilon_2(t) = \sin\left[\frac{\pi}{2} \cdot t\right] - f_2\left[\frac{\pi}{2} \cdot t\right]$$

$$= \left[-10 + 3\pi + \frac{\pi^2}{8} - \frac{\pi^3}{48}\right]t^3 + \left[15 - 4\pi - \frac{\pi^2}{4}\right]t^4 + \left[-6 + \frac{3\pi}{2} + \frac{\pi^2}{8} + \frac{\pi^5}{3840}\right]t^5 - \frac{1}{7!}\left[\frac{\pi}{2} \cdot t\right]^7 + \frac{1}{9!}\left[\frac{\pi}{2} \cdot t\right]^9 + \ldots \quad (92)$$

$$= e_3 t^3 + e_4 t^4 + e_5 t^5 + e_7 t^7 + e_9 t^9 + \ldots$$

where $e_3 = -10 + 3\pi + \pi^2/8 - \pi^3/48 = 0.0125144$ and $e_4 = 15 - 4\pi - \pi^2/4 = -0.0337717$ etc. The approach is then to proceed in a manner consistent with the first order case detailed in Appendix 1 and the following form is considered:

$$\varepsilon_2(t) = c_3 t^3 (1-t)^{p_3} + c_4 t^4 (1-t)^{p_4} + \ldots + c_n t^n (1-t)^{p_n} + \ldots$$

$$= c_3 t^3 \left[1 - p_3 t + \binom{p_3}{2} t^2 - \ldots \pm t^{p_3}\right] + c_4 t^4 [1 - p_4 t + \ldots \pm t^{p_4}] + c_5 t^5 [1 - p_5 t + \ldots \pm t^{p_5}] + \ldots \quad (93)$$

$$= c_3 t^3 + [-p_3 c_3 + c_4] t^4 + \left[\binom{p_3}{2} c_3 - p_4 c_4 + c_5\right] t^5 + \ldots$$

Equating Equation 92 and Equation 93, it immediately follows that $c_3 = e_3 > 0$. For the coefficient $c_4$ to be positive, the factor $p_3$ needs to be chosen such that $e_4 = -p_3 c_3 + c_4$ and $c_4 = e_4 + p_3 c_3 > 0$. This is guaranteed when $p_3 > -e_4/c_3$. Choosing the smallest integer greater than zero such this is the case yields $p_3 = \lfloor -e_4/e_3 \rfloor + 1$ and for the defined values of $e_3$ and $e_4$ the required value is $p_3 = 3$. It then follows that $c_4 = -15 + 5\pi + \pi^2/8 - \pi^3/16 = 0.00377153$.

The ratio of $c_4$ to $c_3$ then is

$$\frac{c_4}{c_3} = \frac{e_4 + p_3 c_3}{e_3} = \frac{e_4}{e_3} + p_3 = \frac{e_4}{e_3} + \left\lfloor\frac{-e_4}{e_3}\right\rfloor + 1 = 1 - \left[\frac{-e_4}{e_3} - \left\lfloor\frac{-e_4}{e_3}\right\rfloor\right] \quad (94)$$

which implies $0 < c_4/c_3 < 1$. Proceeding in this manner yields the stated form and proof of positive coefficients with decreasing values.

The coefficients are defined according to

$$c_3 = -10 + 3\pi + \frac{\pi^2}{8} - \frac{\pi^3}{48} \qquad c_4 = -15 + 5\pi + \frac{\pi^2}{8} - \frac{\pi^3}{16} \quad (95)$$

$$c_5 = -36 + \frac{25\pi}{2} + \frac{\pi^2}{4} - \frac{3\pi^3}{16} + \frac{\pi^5}{3840} \qquad c_6 = -64 + 23\pi + \frac{3\pi^2}{8} - \frac{19\pi^3}{48} + \frac{\pi^5}{960} \quad (96)$$

$$c_7 = -36 + 14\pi + \frac{\pi^2}{8} - \frac{5\pi^3}{16} + \frac{\pi^5}{640} - \frac{\pi^7}{645120} \qquad c_8 = -45 + 18\pi + \frac{\pi^2}{8} - \frac{7\pi^3}{16} + \frac{\pi^5}{384} - \frac{\pi^7}{215040} \quad (97)$$

$$c_9 = -100 + \frac{81\pi}{2} + \frac{\pi^2}{4} - \frac{49\pi^3}{48} + \frac{5\pi^5}{768} - \frac{\pi^7}{71680} + \frac{\pi^9}{185,794,560} \quad (98)$$

$$c_{10} = -166 + 68\pi + \frac{3\pi^2}{8} - \frac{85\pi^3}{48} + \frac{23\pi^5}{1920} - \frac{19\pi^7}{645120} + \frac{\pi^9}{46,448,640} \quad (99)$$

$$c_{11} = -78 + 33\pi + \frac{\pi^2}{8} - \frac{15\pi^3}{16} + \frac{7\pi^5}{960} - \frac{\pi^7}{43008} + \frac{\pi^9}{30,965,760} - \frac{\pi^{11}}{81,749,606,400} \tag{100}$$

$$c_{12} = -91 + 39\pi + \frac{\pi^2}{8} - \frac{55\pi^3}{48} + \frac{3\pi^5}{320} - \frac{\pi^7}{30720} + \frac{\pi^9}{18,579,456} - \frac{\pi^{11}}{27,249,868,800} \tag{101}$$

To find a general formula for the coefficients consider the definitions:

$$c_0 = -1 \qquad c_1 = -2 + \frac{\pi}{2} = 2c_0 + \frac{\pi}{2} \qquad c_2 = -3 + \pi + \frac{\pi^2}{8} = 2c_1 - c_0 + \frac{\pi^2}{8} \tag{102}$$

It is then the case that

$$c_3 = -10 + 3\pi + \frac{\pi^2}{8} - \frac{\pi^3}{48} = c_2 + 4c_1 - c_0 - \frac{\pi^3}{2^3 3!} \tag{103}$$

and it can then be readily shown that

$$c_4 = 3c_3 - 2c_2 - 4c_1 - c_0 \qquad c_5 = 3c_4 - c_2 - 3c_1 + \frac{\pi^5}{2^5 5!} \tag{104}$$

$$c_6 = 4c_5 - 6c_4 + c_3 \qquad c_7 = 2c_6 - 2c_5 - 2c_4 + c_3 - \frac{\pi^7}{2^7 7!} \tag{105}$$

$$c_8 = 3c_7 - 3c_6 + 4c_5 - c_4 \qquad c_9 = 4c_8 - 3c_7 + c_6 - c_5 + \frac{\pi^9}{2^9 9!} \tag{106}$$

The iterative formulas then follow:

$$\begin{aligned} c_k &= 4c_{k-1} - 6c_{k-2} + c_{k-3} & k \in \{6, 10, 14, 18, \ldots\} \\ c_k &= 2c_{k-1} - 2c_{k-2} - 2c_{k-3} + c_{k-4} - \frac{\pi^k}{2^k k!} & k \in \{7, 11, 15, 19, \ldots\} \\ c_k &= 3c_{k-1} - 3c_{k-2} + 4c_{k-3} - c_{k-4} & k \in \{8, 12, 16, 20, \ldots\} \\ c_k &= 4c_{k-1} - 3c_{k-2} + c_{k-3} - c_{k-4} + \frac{\pi^k}{2^k k!} & k \in \{9, 13, 17, 21, \ldots\} \end{aligned} \tag{107}$$

The proof that the series defining $\varepsilon_2$ converges, and is bounded below by zero, follows in an identical manner to that given in Appendix 1.

## Appendix 3:   Proof of Theorem 5.2

Equation 29 implies

$$\sin(x) = f_1(x) + \varepsilon_1\left[\frac{2x}{\pi}\right] \qquad x \in \left[0, \frac{\pi}{2}\right] \tag{108}$$

where $f_1$ is defined by Equation 24 and the convergent series for $\varepsilon_1$ is defined in Theorem 4.1. As $f_1$ can be written in the form

$$f_1(x) = \frac{2x}{\pi} + \frac{2x}{\pi}\left[1 - \frac{2x}{\pi}\right] + c_1\left[\frac{2x}{\pi}\right]\left[1 - \frac{2x}{\pi}\right]^2 \qquad c_1 = 2\left[-1 + \frac{\pi}{4}\right] \tag{109}$$

the required result follows, namely

$$\sin(x) = \frac{2x}{\pi} + \frac{2x}{\pi}\left[1 - \frac{2x}{\pi}\right] + \sum_{k=1}^{\infty} c_k\left[\frac{2x}{\pi}\right]^k\left[1 - \frac{2x}{\pi}\right]^p \qquad p = \begin{cases} 2 & k = 1, 2, 5, 6, 9, \ldots \\ 3 & k = 3, 4, 7, 8, \ldots \end{cases} \tag{110}$$

## Appendix 4: Proof of Theorem 5.3

Equation 34 implies

$$\sin(x) = f_2(x) + \varepsilon_2\left[\frac{2x}{\pi}\right] \qquad x \in \left[0, \frac{\pi}{2}\right] \tag{111}$$

where $f_2$ is defined by Equation 25 and the convergent series for $\varepsilon_2$ is defined in Theorem 4.2, As $f_2$ can be written in the form

$$f_2(x) = 1 - \frac{\pi^2}{8}\left[1 - \frac{2x}{\pi}\right]^2 + \left[-1 + \frac{\pi^2}{8}\right]\left[1 - \frac{2x}{\pi}\right]^4 + \left[-4 + \frac{\pi}{2} + \frac{\pi^2}{4}\right]\left[\frac{2x}{\pi}\right]\left[1 - \frac{2x}{\pi}\right]^4 + \\ \left[-10 + 2\pi + \frac{3\pi^2}{8}\right]\left[\frac{2x}{\pi}\right]^2\left[1 - \frac{2x}{\pi}\right]^3 \tag{112}$$

the required result follows, namely:

$$\sin(x) = 1 - \frac{\pi^2}{8}\left[1 - \frac{2x}{\pi}\right]^2 + \sum_{k=0}^{\infty} c_k\left[\frac{2x}{\pi}\right]^k\left[1 - \frac{2x}{\pi}\right]^p$$

$$p = \begin{cases} 3 & k = 2, 3, 6, 7, 10, 11, 14, 15, , \ldots \\ 4 & k = 0, 1, 4, 5, 8, 9, 12, 13, 16, 17, \ldots \end{cases} \tag{113}$$

$$c_0 = -1 + \frac{\pi^2}{8} \qquad c_1 = -4 + \frac{\pi}{2} + \frac{\pi^2}{4} \qquad c_2 = -10 + 2\pi + \frac{3\pi^2}{8}$$

## Appendix 5: Proof of Theorem 5.4

First, the definitions $f_1$ and $f_2$, as defined by Equation 24 and Equation 25, imply

$$2f_2(x) - f_1(x) = 2\left[x + \frac{80}{\pi^3}\left[1 - \frac{3\pi}{10} - \frac{\pi^2}{80}\right]x^3 - \frac{240}{\pi^4}\left[1 - \frac{4\pi}{15} - \frac{\pi^2}{60}\right]x^4 + \frac{192}{\pi^5}\left[1 - \frac{\pi}{4} - \frac{\pi^2}{48}\right]x^5\right] - \\ \left[x + \frac{12}{\pi^2}\left[1 - \frac{\pi}{3}\right]x^2 - \frac{16}{\pi^3}\left[1 - \frac{\pi}{4}\right]x^3\right] \tag{114}$$

$$= x - \frac{12}{\pi^2}\left[1 - \frac{\pi}{3}\right]x^2 + \frac{176}{\pi^3}\left[1 - \frac{13\pi}{44} - \frac{\pi^2}{88}\right]x^3 + \frac{480}{\pi^4}\left[1 - \frac{4\pi}{15} - \frac{\pi^2}{60}\right]x^4 + \frac{384}{\pi^5}\left[1 - \frac{\pi}{4} - \frac{\pi^2}{48}\right]x^5$$

Second, utilizing the transformed functions associated with the interval $(0, 1)$ rather than the interval $(0, \pi/2)$, the requirement is for $\varepsilon_1(t) > 2\varepsilon_2(t)$, $t \in (0, 1)$, i.e. for

$$c_2 t^2 (1-t)^2 + c_3 t^3 (1-t)^3 + c_4 t^4 (1-t)^3 + c_5 t^5 (1-t)^2 + c_6 t^6 (1-t)^2 + c_7 t^7 (1-t)^3 + \ldots \\ > 2[d_3 t^3 (1-t)^3 + d_4 t^4 (1-t)^4 + d_5 t^5 (1-t)^4 + d_6 t^6 (1-t)^3 + d_7 t^7 (1-t)^3 + d_8 t^8 (1-t)^4 + \ldots] \tag{115}$$

where the positive coefficients $c_k$ and $d_k$ are defined, respectively, in Theorem 4.1 and Theorem 4.2. Consider

$$c_k t^k (1-t)^{p_1} - 2d_{k+1} t^{k+1} (1-t)^{p_1+1} \qquad p_1 \in \{2, 3\}, t \in (0, 1), k \in \{2, 3, 4, \ldots\}$$

$$= [c_k - 2d_{k+1} t(1-t)] t^k (1-t)^{p_1} = [c_k - 2d_{k+1}] t^k (1-t)^{p_1} + 2d_{k+1}[1 - t(1-t)] t^k (1-t)^{p_1} \tag{116}$$

$$> [c_k - 2d_{k+1}] t^k (1-t)^{p_1}$$

Hence, a sufficient condition for $\varepsilon_1(t) > 2\varepsilon_2(t)$, $t \in (0, 1)$, is for $c_k - 2d_{k+1} > 0$, $k \in \{2, 3, 4, \ldots\}$. This is easy to confirm, along with the rapid decline toward zero in the sequence $c_k - 2d_{k+1}$, $k \in \{2, 3, 4, \ldots\}$.